\documentclass[12pt,twoside,reqno]{amsart}
\usepackage[colorlinks=true,citecolor=blue]{hyperref}
\usepackage{mathptmx, amsmath, amssymb, amsfonts, amsthm, mathptmx, enumerate, color}
\usepackage{graphicx}
\usepackage{epstopdf}

\newtheorem{theorem}{Theorem}[section]

\theoremstyle{definition}

\numberwithin{equation}{section}

\usepackage{tikz}

\DeclareMathOperator*{\st}{s.t.}

\usepackage{pgf}
\usepackage{pgfplots}
\pgfplotsset{compat = newest}
\usepackage{caption}

\usepackage{algorithm}
\usepackage{algpseudocode}
\newcounter{protocol}
\makeatletter

\DeclareMathOperator{\E}{\mathbb{E}}

\begin{document}
\setcounter{page}{1}

\vspace*{1.0cm}
\title[Convergence analysis of coordinate descent]
{Convergence rate analysis of randomized and cyclic coordinate descent for convex optimization  through semidefinite programming}
\author[H. Abbaszadehpeivasti, E. de Klerk, M. Zamani]{ Hadi Abbaszadehpeivasti$^{1}$, Etienne de Klerk$^{1,*}$, Moslem Zamani$^1$}
\maketitle
\vspace*{-0.6cm}

\begin{center}
{\footnotesize {\it

$^1$Tilburg University, Department of Econometrics and Operations Research, Tilburg, The Netherlands

}}\end{center}

\vskip 4mm {\small\noindent {\bf Abstract.}
In this paper, we study randomized and cyclic coordinate descent for convex unconstrained optimization problems.
 We improve the known convergence rates in some cases  by using the numerical semidefinite programming performance estimation method.
 As a spin-off we provide a method to analyse the worst-case performance of the Gauss-Seidel iterative method for linear systems where
 the coefficient matrix is positive semidefinite with a positive diagonal.


}

\renewcommand{\thefootnote}{}
\footnotetext{ $^*$Corresponding author.
\par
E-mail addresses: h.abbaszadehpeivasti@tilburguniversity.edu (H. Abbaszadehpeivasti), e.deklerk@tilburguniversity.edu (E. de Klerk), m.zamani\_1@tilburguniversity.edu (M. Zamani).
\par
Received ...; Accepted February .... }

\section{Introduction}

We consider the unconstrained optimization problem
\begin{align}\label{P}
f^{\star}=\min_{x\in \mathbb{R}^n} f(x),
\end{align}
where $f: \mathbb{R}^n\to \mathbb{R}$ is convex. We assume that $f$ attains its minimum and $f^\star$ denotes the optimal value. In addition, we assume that $f$ is an $L$-smooth function, that is,
$$
\|\nabla f(y)-\nabla f(x)\|\leq L\|y-x\|, \ \ \ \forall y, x\in\mathbb{R}^n.
$$
Moreover, we denote the component Lipschitz constants by $\ell_i$ ($i \in \{1,\ldots,n\}$), i.e.,
\begin{align}\label{Co.Lip}
  |[\nabla f(x+te_i)]_i-[\nabla f(x)]_i|\leq \ell_i|t|, \ \ \ \forall x\in\mathbb{R}^n, t\in\mathbb{R},
\end{align}
where $e_i$ is the $i$th standard unit vector. Let $\ell_{\max}:=\max_{1\leq i\leq n} \ell_i$, and note that $1\leq \tfrac{L}{\ell_{\max}}\leq n$.

Due to the  simplicity and small per-iteration cost,    coordinate descent methods have been employed extensively for large-scale optimization problems
\cite{nesterov2012efficiency, wright2015coordinate}.

The generic coordinate descent method is shown in Algorithm \ref{Alg0}.

\begin{algorithm}
\caption{Generic coordinate descent}
\begin{algorithmic}
\State Set $N$ and $\{t_k\}_{k=0}^{N-1}$  (step lengths) and pick $x^0\in\mathbb{R}^n$.
\State For $k=0,1, \ldots, N-1$ perform the following step:\\
\begin{enumerate}
\item
Choose an index $i_k$ from $\{1, 2, . . . , n\}$.
\item
$x^{k+1}=x^k-t_k[\nabla f(x^k)]_{i_k} e_{i_k}$.
\end{enumerate}
\end{algorithmic}
\label{Alg0}
\end{algorithm}

 In this paper, we revisit the worst-case convergence rate analysis for
Algorithm \ref{Alg0} for two of the best known variants, namely \textit{randomized coordinate descent}, and \textit{cyclic coordinate descent}.
In the former, the index $i_k$ is chosen uniformly at random from $\{1, 2, . . . , n\}$, and in the latter, the cyclic ordering is used.

We will improve the best-known convergence rates from the literature for some specific values of
the parameters $n,L,N, t_k$ for $k\in\{0,1,\ldots,N-1\}$ and $\ell_i$ for $i\in\{1,\ldots,n\}$.
Finally, the Gauss-Seidel iterative method for positive semidefinite linear systems is a special case cyclic coordinate descent for convex quadratic
functions, and we will investigate the implications of our analysis for this classical method as well.

Recently, Kamri et al.\ \cite{kamri2022worst} studied the convergence of the coordinate descent algorithm using the semidefinite programming (SDP) performance estimation method, that was introduced by Drori and Teboulle \cite{drori2014performance}.
We will also use SDP performance estimation in our analysis, and our main contribution may be seen as the extension and refinement of the approach by Kamri et al.\ \cite{kamri2022worst}.
SDP performance estimation has been applied to the analysis of many iterative methods (other than coordinate descent); the interested reader may
consult \cite{taylor2017smooth,de2020worst,drori2020complexity,abbaszadehpeivasti2021rate,TeodorHypo}  and the references therein. For general background information on SDP, see e.g.\ \cite{wolkowicz2012handbook}.

\subsection*{Notation and background results}

 We use $\langle \cdot, \cdot\rangle$ and $\| \cdot\|$ to denote the Euclidean inner product and norm, respectively, unless indicated otherwise. The column vector $e_i$ represents the $i$-th standard unit vector and $I$ stands for the identity matrix.
  For a matrix $A$, $a_{ij}$ denotes its $(i, j)$-th entry,
  and $A^\top$ represents the transpose of $A$.
The function $f:\mathbb{R}^n\to\mathbb{R}$ is called $\mu$-strongly convex function if the function $x \mapsto f(x)-\tfrac{\mu}{2}\| x\|^2$ is convex.
 Clearly, any convex function is $0$-strongly convex.
  We denote the set of real-valued convex functions which are $L$-smooth and $\mu$-strongly convex by $\mathcal{F}_{\mu,L}(\mathbb{R}^n)$.

Let $\mathcal{I}$ be a finite index set and let $\{x^i; g^i; f^i\}_{i\in \mathcal{I}}\subseteq \mathbb{R}^n\times \mathbb{R}^n\times \mathbb{R}$.
A set $\{x^i; g^i; f^i\}_{i\in \mathcal{I}}$ is called $\mathcal{F}_{\mu,L}$-interpolable if there exists $f\in\mathcal{F}_{\mu,L}(\mathbb{R}^n)$
 with
$$
f(x^i)=f^i, \ g^i\in\partial f(x^i) \ \ i\in\mathcal{I}.
$$
The next theorem gives necessary and sufficient conditions for $\mathcal{F}_{\mu,L}$-interpolablity.
\begin{theorem}\cite{taylor2017smooth}\label{T1}
Let $L\in (0, \infty)$ and $\mu\in [0, \infty)$ and $f\in\mathcal{F}_{\mu,L}(\mathbb{R}^n)$. For any $x,y\in\mathbb{R}^n$, we have
{\small{
\begin{align}\label{inter.1}
\tfrac{1}{2(1-\tfrac{\mu}{L})}\left(\tfrac{1}{L}\left\|\nabla f(x)-\nabla f(y)\right\|^2+\mu\left\|x-y\right\|^2-\tfrac{2\mu}{L}\left\langle \nabla f(y)-\nabla f(x),y-x\right\rangle\right)\leq f(x)-f(y)-\left\langle \nabla f(y), x-y\right\rangle.
\end{align}
}}
Conversely, if $\mathcal{I}$ is a finite index set and $\{x^i; g^i; f^i\}_{i\in \mathcal{I}}\subseteq \mathbb{R}^n\times \mathbb{R}^n\times \mathbb{R}$ are given data, then the data set is $\mathcal{F}_{\mu,L}$-interpolable if it satisfies \eqref{inter.1} in the sense that, for each pair $i,j \in \mathcal{I}$:
\[
\tfrac{1}{2(1-\tfrac{\mu}{L})}\left(\tfrac{1}{L}\left\|g^i-g^j\right\|^2+\mu\left\|x^i-x^j\right\|^2-\tfrac{2\mu}{L}\left\langle g^j-g^i,x^j-x^i\right\rangle\right)\leq f^i-f^j-\left\langle g^j, x^i-x^j\right\rangle.
\]

\end{theorem}

\section{Convergence rate of randomized coordinate descent}
The randomized coordinate descent method is shown in Algorithm \ref{Alg1} for easy reference.

\begin{algorithm}
\caption{Randomized coordinate descent}
\begin{algorithmic}
\State Set $N$ and $\{t_k\}_{k=0}^{N-1}$ (step lengths) and pick $x^0\in\mathbb{R}^n$.
\State For $k=0,1, \ldots, N-1$ perform the following step:\\
\begin{enumerate}
\item
Choose index $i_k$ with uniform probability from $\{1, 2, . . . , n\}$.
\item
$x^{k+1}=x^k-t_k[\nabla f(x^k)]_{i_k} e_{i_k}$.
\end{enumerate}
\end{algorithmic}
\label{Alg1}
\end{algorithm}

We proceed to revisit its worst-case convergence rate for three classes for function, namely convex $L$-smooth functions, convex quadratic functions, and strongly convex, $L$-smooth functions.

\subsection{The case of $L$-smooth functions}
Regarding the convergence of Algorithm \ref{Alg1} for $L$-smooth convex functions, the following is known. (We state the result as in the survey \cite[Theorem 1]{wright2015coordinate}, but it is originally due to Nesterov \cite{nesterov2012efficiency}).
\begin{theorem}\cite[Theorem 1]{wright2015coordinate}
Let $f:\mathbb{R}^n\to\mathbb{R}$ is an $L$-smooth convex function for some $L>0$. If $t_k=\tfrac{1}{\ell_{\max}}$ for all $k$, then, for each $k > 0$,
\begin{align}\label{Wrigt}
  \mathbb{E}\left( f(x^{k}) \right)-f^\star \leq \left(\tfrac{2n\ell_{\max}}{k}\right) R_0^2,
\end{align}
where $R_0$ satisfies $\max_{x^{\star}\in\mathbb{S}}\max_{x}\{\|x-x^{\star}\|:f(x)\leq f(x^0)\}\leq R_0$ and  $\mathbb{S}$ denotes
the optimal solution set.
\end{theorem}

 In this section, we study the behaviour of randomized coordinate descent method for $L$-smooth convex functions. The worst-case convergence rate of the Algorithm \ref{Alg1} can be formulated as follows.
\begin{align}\label{P1}
\nonumber   \max & \ \mathbb{E}[f(x^N)] -f(x^\star)\\
 \nonumber \st   &  \  f\in\mathcal{F}_{0,L}(\mathbb{R}^n) \\
  & \text{$f$ satisfies \eqref{Co.Lip} for  every $x\in\mathbb{R}^n$ for some $\ell_i$, $i\in\{1,\dots,n\}$}\\
   \nonumber & \|x^0-x^\star\|^2\leq\Delta\\
  \nonumber   & \ \text{$x^k$ $k \in \{1,2,...,N\} $ are generated by Algorithm \ref{Alg1} with respect to $x^0$ and step length $t_k$}\\
\nonumber &  x^0\in\mathbb{R}^n, \ \nabla f(x^\star)=0,
\end{align}
where $f, x^k, x^\star$ are decision variables and $t, L, n$ and $\ell_i, i\in\{1,\dots,n\}$, are the given parameters. Problem \eqref{P1} in general is intractable. Moreover, note that $x^k$ depends on the index $i_k$ which is chosen uniformly at random from the set $\{1,\dots,n\}$ therefore \eqref{P1} is a stochastic programming problem. To deal with this we introduce a random variable $d^k$ which depends on the index $i_k$ and is defined by $d^k :=[\nabla f(x^k)]_{i_k} e_{i_k}$. Note that $d^k$ has the following properties:
\begin{align}\label{d.P}
   \nonumber & \E\left[\|d^k\|^2\right]=\tfrac{1}{n}\E\left[\|\nabla f(x^k)\|^2\right]\\
   & \E\left[\left\langle d^k,\nabla f(x^k)\right\rangle\right]=\tfrac{1}{n}\E\left[\|\nabla f(x^k)\|^2\right]\\
   \nonumber & \E\left[\left\langle d^k,x^k\right\rangle\right] =\tfrac{1}{n}\E\left[\left\langle \nabla f(x^k),x^k\right\rangle\right],
\end{align}
where the expectation again refers to the joint distribution of all the random variables $d^k$ for $k\in\{0,1,\ldots,N\}$ and $x^k, \nabla f(x^k), f(x^k)$ for $k\in\{0,1,\ldots,N\}$.
By  Taylor's theorem and \eqref{Co.Lip}, we have
\begin{align}\label{re.co.lip}
   & f(x^{k+1})\leq f(x^k)+\left\langle \nabla f(x^{k}),x^{k+1}-x^k\right\rangle+ \tfrac{\ell_{\max}}{2}\|x^{k+1}-x^k\|^2 \\
\nonumber   & f(x^{k})\leq f(x^{k+1})+\left\langle \nabla f(x^{k+1}),x^{k}-x^{k+1}\right\rangle+ \tfrac{\ell_{\max}}{2}\|x^k-x^{k+1}\|^2,
\end{align}
where, $\ell_{\max}=\max_{i\in\{1,\dots,n\}} \ell_i$ as before. Therefore, the relaxation of problem \eqref{P1} is given by
\begin{align}\label{P2}
\nonumber   \max & \  \mathbb{E}[f(x^N)] -f(x^\star)\\
 \nonumber \st   &  \  \{(x^k; \nabla f(x^k); f(x^k))\} \ \textrm{satisfy \eqref{inter.1} for $k\in\{0,1,\cdots,N, \star\}$ w.r.t.}\ \mu=0, L \\
  \nonumber    &  \  \{(x^k; \nabla f(x^k); f(x^k))\} \ \textrm{satisfy \eqref{re.co.lip} for $k\in\{0,1,\cdots,N\}$ w.r.t.}\ \ell_{\max} \\
  & \ \{x^k;\nabla f(x^k); d^k\} \ \textrm{satisfies \eqref{d.P}} \;\;\;  (k\in\{0,\cdots,N\})\\
\nonumber    & \|x^0-x^\star\|^2\leq\Delta\\
  \nonumber   & \ x^{k+1}=x^{k}-t_kd^k\\
\nonumber &  x^0\in\mathbb{R}^n, \ \nabla f(x^\star)=0,
\end{align}
where $f(x^k), x^k, x^\star, \nabla f(x^k)$ and $d^k$ are decision variables. Note that because the problem \eqref{P} is invariant under translation, without loss of generality we may assume that $x^\star$ is the zero vector. Since $x^{k+1}=x^{k}-t_kd^k$ is a recursive relation, $x^k$ can be written as linear combination of $x^0$ and $d^i$s. In this way, all the unknowns appear as entries in the following matrix:
{\footnotesize{
\begin{align*}
  G=\begin{pmatrix}
      \E[\|x^0\|^2] & \E[\langle x^0,\nabla f(x^0)\rangle] & \cdots & \E[\langle x^0,\nabla f(x^N)\rangle] & \E[\langle x^0, d^0\rangle] & \cdots & \E[\langle x^0, d^N\rangle] \\
      \E[\langle \nabla f(x^0), x^0\rangle] & \E[\|\nabla f(x^0)\|^2] & \cdots & \E[\langle \nabla f(x^0),\nabla f(x^N)\rangle] & \E[\langle \nabla f(x^0), d^0\rangle] & \cdots & \E\langle \nabla f(x^0), d^N\rangle \\
      \vdots &  \vdots & \ddots & \vdots & \vdots & \ddots & \vdots\\
      \E[\langle d^N, x^0\rangle] & \E[\langle d^N, \nabla f(x^0)\rangle] & \cdots & \E[\langle d^N,\nabla f(x^N)\rangle] & \E[\langle d^N, d^0\rangle] & \cdots & \E[\|d^N\|^2] \\
    \end{pmatrix}.
\end{align*}
}}
Note that $G$ is the expectation of a random Gram matrix. Since every realization of this random matrix is positive semidefinite, and the expectation preserves positive semidefiniteness, it follows that $G$ is positive semidefinite as well. Therefore, problem \eqref{P2} can be written as an SDP problem, where the variables are $G$ and $\E [f(x^i)]$.

In what follows we compare the convergence rate derived by solving the problem \eqref{P2} and the bound  by Wright \eqref{Wrigt} for some specific values of the parameters $n,L,N, t_k$ for $k\in\{0,1,\cdots,N\}$ and $\ell_i$ for $i\in\{1,\cdots,n\}$.
All the figures in this paper were obtained by solving the SDP problems with the solver Mosek \cite{mosek}, using the  Yalmip \cite{Lofberg2004} Matlab interface.

\begin{figure}[H]
\begin{tikzpicture}
\begin{axis}[
    xmin = 1, xmax = 30,
    ymin = 0, ymax = 4,
    xtick distance = 5,
    ytick distance = 2,
    grid = both,
    minor tick num = 2,
    major grid style = {lightgray},
    minor grid style = {lightgray!50},
    width = \textwidth,
    height = 0.5\textwidth,
    xlabel = {Number of iterations $(N)$},
    ylabel = {$\mathbb{E}[f(x^N)] -f(x^\star)$},]

\addplot[
    domain = 1:30,
    samples = 1000,
    smooth,
    thick,
    blue,
] {(4)/(x)};
\addplot[
    domain = 1:30,
    samples = 1000,
    smooth,
    thick,
    red,
] {(8)/(x)};
\addplot[
    smooth,
    thin,
    blue,
    dashed
] table[col sep=comma]  {fign2L2.dat};
\addplot[
    smooth,
    thin,
    red,
    dashed
] table[col sep=comma]  {fign4L2.dat};

 \legend{
    Bound \eqref{Wrigt} ($n=2$),
    Bound \eqref{Wrigt} ($n=4$),
    PEP bound \eqref{P2} ($n=2$),
    PEP bound \eqref{P2} ($n=4$),
}

\end{axis}

\end{tikzpicture}
\caption{Convergence rate for Algorithm \ref{Alg1} computed by performance estimation problem \eqref{P2} (dashed lines)
and the bound given by \eqref{Wrigt} (thick lines)
for $L=2, l_{\max}=1, t=1, \Delta=1$ and different $n$}
\label{fig1}
\end{figure}

\begin{figure}[H]
\begin{tikzpicture}
\begin{axis}[
    xmin = 1, xmax = 30,
    ymin = 0, ymax = 4,
    xtick distance = 5,
    ytick distance = 2,
    grid = both,
    minor tick num = 2,
    major grid style = {lightgray},
    minor grid style = {lightgray!50},
    width = \textwidth,
    height = 0.5\textwidth,
    xlabel = {Number of iterations $(N)$},
    ylabel = {$\mathbb{E}[f(x^N)] -f(x^\star)$},]

\addplot[
    domain = 1:30,
    samples = 1000,
    smooth,
    thick,
    blue,
] {(8)/(x)};
\addplot[
    domain = 1:30,
    samples = 1000,
    smooth,
    thick,
    red,
] {(16)/(x)};
\addplot[
    smooth,
    thin,
    blue,
    dashed
] table[col sep=comma]  {fign2L4.dat};
\addplot[
    smooth,
    thin,
    red,
    dashed
] table[col sep=comma]  {fign4L4.dat};

 \legend{
    Bound \eqref{Wrigt} ($n=2$),
    Bound \eqref{Wrigt} ($n=4$),
    PEP bound \eqref{P2} ($n=2$),
    PEP bound \eqref{P2} ($n=4$),
}

\end{axis}

\end{tikzpicture}
\caption{Convergence rate for Algorithm \ref{Alg1} computed by performance estimation problem \eqref{P2} (dashed lines)
and the bound given by \eqref{Wrigt} (thick lines)
for $L=4, l_{\max}=2, t=0.5, \Delta=1$ and different $n$}
\label{fig1}
\end{figure}
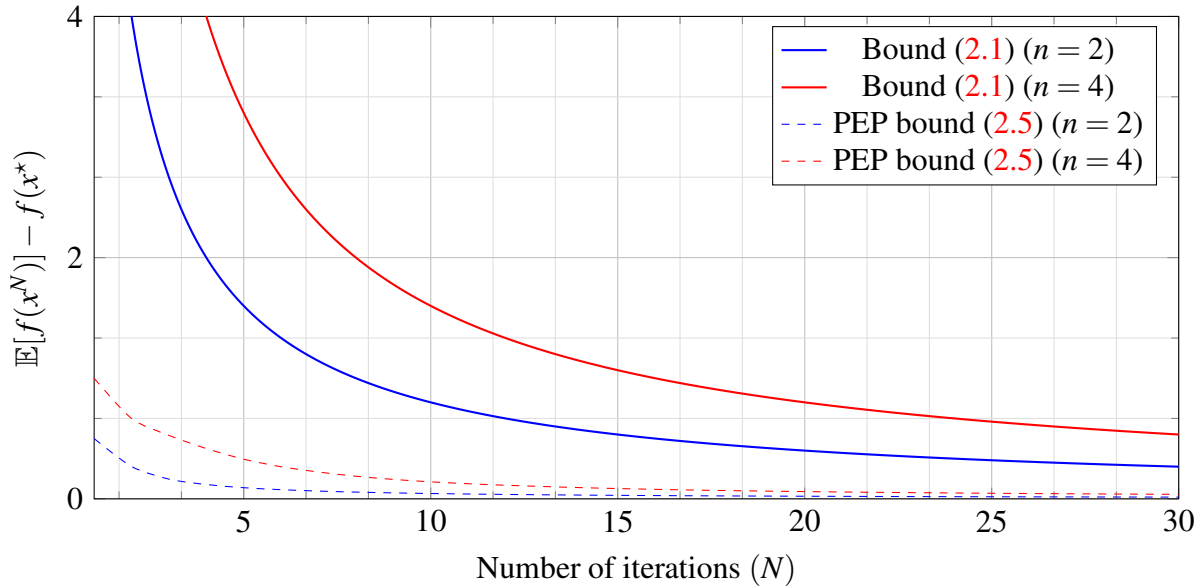
Note that the convergence rate provided by solving performance estimation is strictly better than the bound given by Wright. In other words, the bound \eqref{Wrigt} is not tight for the values of the parameters that we considered. Moreover, the bound given by performance estimation can also be calculated for different step lengths than the fixed  step lengths $1/\ell_{\max}$ in the bound \eqref{Wrigt}.

\subsection{The case of convex quadratic functions}

In this section we study the convergence of rate of the randomized coordinate descent method in case that the objective function is quadratic function of the form
\begin{align}\label{QP}
\min_{x\in \mathbb{R}^n} f(x):=\tfrac{1}{2}x^\top Ax-b^\top x,
\end{align}
where, $A$ is a symmetric positive semidefinite matrix. To study this case we need to add additional constraint to restrict our model to quadratic functions. The following necessary condition for $f$ to be a quadratic function can be verified easily, and has been used in SDP performance analysis by Drori et al \cite{drori2020complexity}:
\begin{align*}
  \tfrac{1}{2}\langle \nabla f(x)-\nabla f(y),x-y\rangle=f(x)-f(y)-\langle \nabla f(y),x-y\rangle.
\end{align*}
Since this constraint holds for every point in the domain we just consider the relaxed constraint that only holds for the point generated by the method in addition to the initial point and the optimal point. In this case we add the following constraint to the problem \eqref{P2}.
\begin{align}\label{qua_cons}
  \tfrac{1}{2}\langle \nabla f(x^i)-\nabla f(x^j),x^i-x^j\rangle=f(x^i)-f(x^j)-\langle \nabla f(x^j),x^i-x^j\rangle \ \ \ \forall i,j\in\{0,\ldots,N,*\}.
\end{align}
In what follows we compare the convergence rate of the randomized coordinate descent method for the general problem \eqref{P2} to the convergence rate for the quadratic problems.

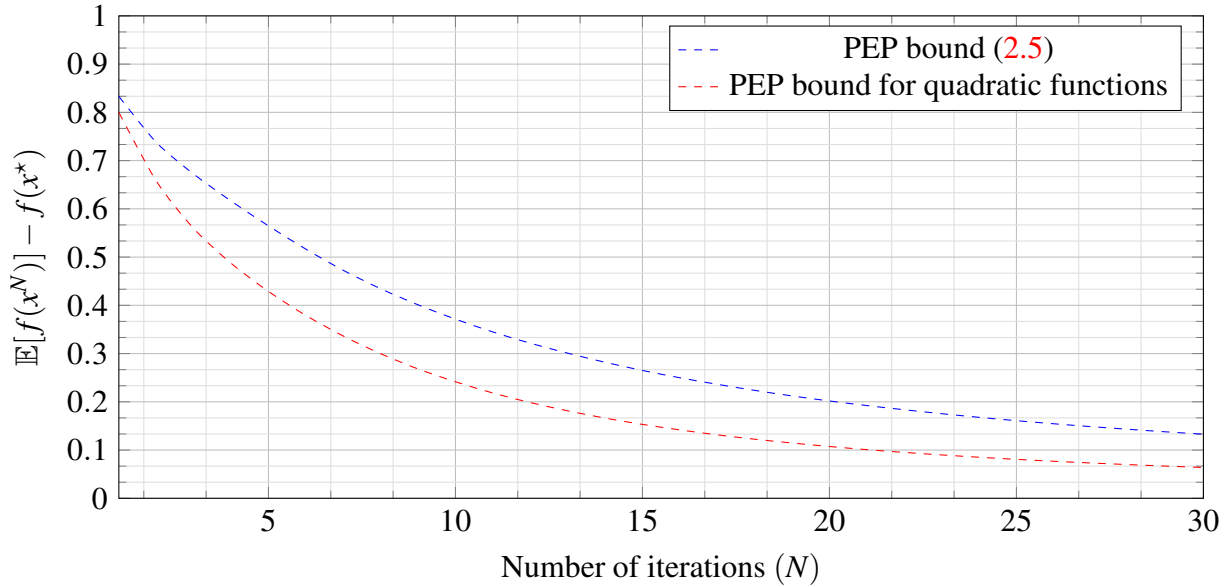
\begin{figure}[H]
\begin{tikzpicture}
\begin{axis}[
    xmin = 1, xmax = 30,
    ymin = 0, ymax = 1,
    xtick distance = 5,
    ytick distance = .1,
    grid = both,
    minor tick num = 2,
    major grid style = {lightgray},
    minor grid style = {lightgray!50},
    width = \textwidth,
    height = 0.5\textwidth,
    xlabel = {Number of iterations $(N)$},
    ylabel = {$\E[f(x^N)] -f(x^\star)$},]

\addplot[
    smooth,
    thin,
    blue,
    dashed
] table  {fig_n10.dat};
\addplot[
    smooth,
    thin,
    red,
    dashed
] table  {figq_n10.dat};

 \legend{
    PEP bound \eqref{P2},
    PEP bound for quadratic functions,
}

\end{axis}

\end{tikzpicture}
\caption{Convergence rate for Algorithm \ref{Alg1} computed by performance estimation problem for quadratic functions (red line)
and the bound given by \eqref{P2} (blue line)
for $n=10, L=2, l_{\max}=1, t=1, \Delta=1$.}
\label{fig3}
\end{figure}

\begin{figure}[H]
\begin{tikzpicture}
\begin{axis}[
    xmin = 1, xmax = 30,
    ymin = 0, ymax = .5,
    xtick distance = 5,
    ytick distance = .1,
    grid = both,
    minor tick num = 2,
    major grid style = {lightgray},
    minor grid style = {lightgray!50},
    width = \textwidth,
    height = 0.5\textwidth,
    xlabel = {Number of iterations $(N)$},
    ylabel = {$\E[f(x^N)] -f(x^\star)$},]

\addplot[
    smooth,
    thin,
    blue,
    dashed
] table  {fig_n2.dat};
\addplot[
    smooth,
    thin,
    red,
    dashed
] table  {figq_n2.dat};

 \legend{
    PEP bound \eqref{P2},
    PEP bound for quadratic functions,
}

\end{axis}

\end{tikzpicture}
\caption{Convergence rate for Algorithm \ref{Alg1} computed by performance estimation problem for quadratic functions (red line)
and the bound given by \eqref{P2} (blue line)
for $n=2, L=2, l_{\max}=1, t=1, \Delta=1$.}
\label{fig3}
\end{figure}
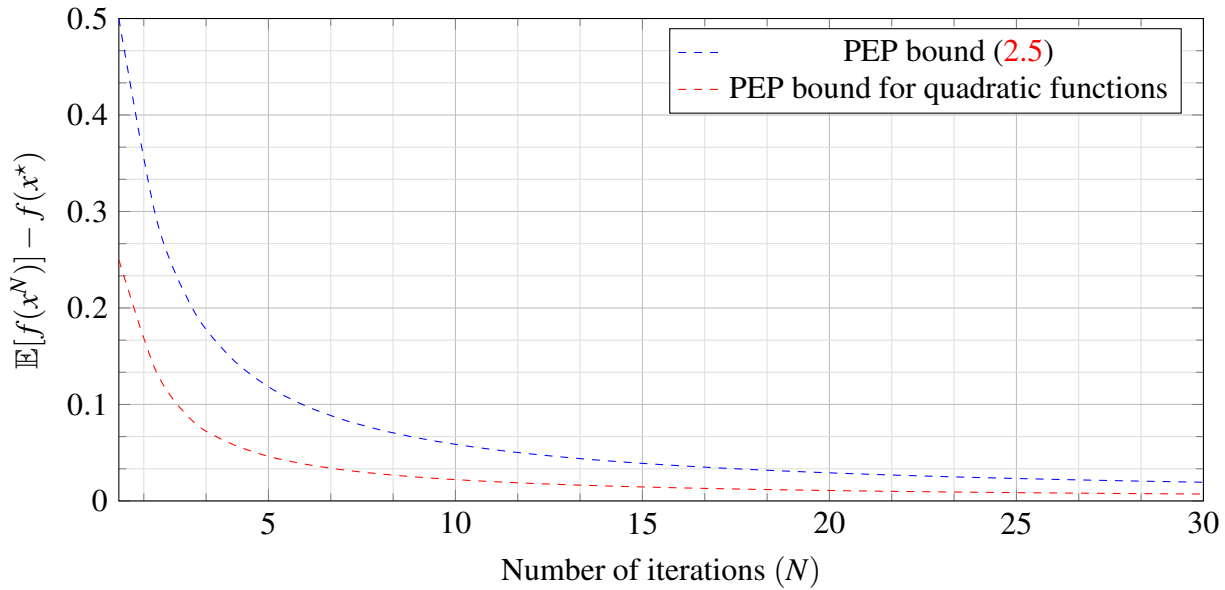
Note that the convergence rate for the quadratic problem is sightly better than that of the general case.

\subsection{The case of $\mu$-strongly convex $L$-smooth functions}
In this section, we study the convergence rate of the $\mu$-strongly convex $L$-smooth functions.
If $\mu>0$ the optimal value of  problem \eqref{P2}  for one iteration of Algorithm \ref{Alg1}, i.e.\ $N=1$,
appears to be the same as the following bound \eqref{SWB} given by Wright \cite{wright2015coordinate}.
\begin{theorem}\cite[Theorem 1]{wright2015coordinate}
Let $f\in\mathcal{F}_{\mu,L}(\mathbb{R}^n)$. If $t_k=\tfrac{1}{\ell_{\max}}$ for each $k$, and $\mu>0$, then, for all $N>0$,
\begin{align}\label{SWB}
  \E\left[f(x^{N}) \right]-f^\star \leq \left(1-\tfrac{\mu}{n\ell_{\max}}\right)^N \left(f(x^{0})-f^\star\right) .
\end{align}
\end{theorem}
The fact that the two bounds seem to coincide does not suggest that  \eqref{SWB} is tight, since the SDP bound \eqref{P2} is a relaxation, and not exact.

\section{Cyclic coordinate descent}
Cyclic coordinate descent is one of the most important coordinate descent algorithms due to its simplicity. The convergence rate of cyclic coordinate descent method for the class of $L$-smooth convex functions is studied by Kamri et al using the performance estimation method \cite{kamri2022worst}. This method is described in Algorithm \ref{Alg-cyc}.
\begin{algorithm}[H]
\caption{Cyclic coordinate descent}
\begin{algorithmic}
\State Set number of cycles $K$, $\{t_k\}_{k=0}^{N-1}$  (step lengths), pick $x^0\in\mathbb{R}^n$ and set $N=nK$.
\State For $k=0,1, 2, \ldots, N-1$ perform the following step:\\
\begin{enumerate}
\item
Set $i=k \pmod n +1$
\item
$x^{k+1}=x^k-t_k[\nabla f(x^k)]_{i} e_{i}$.
\end{enumerate}
\end{algorithmic}
\label{Alg-cyc}
\end{algorithm}
In each iteration the method updates the current point over one of the coordinates in cyclic order.

The following result is known about the rate of convergence. We present it as in \cite{wright2015coordinate}, but it is originally due to Beck and Tetruashvili \cite{beck2013convergence}.
\begin{theorem}\cite[Theorem 3]{wright2015coordinate}
Let $f:\mathbb{R}^n\to\mathbb{R}$ is an $L$-smooth convex function for some $L>0$. If $t_k=\tfrac{1}{\ell_{\max}}$ for all $k$, then, for  $k=n,2n,3n,\ldots$,
\begin{align}\label{Wrigt2}
 \left( f(x^{k}) \right)-f^\star \leq \left(\frac{4n R_0^2\ell_{\max}(1+nL^2/\ell^2_{\max})}{k+8}\right),
\end{align}
where $R_0$ satisfies $\max_{x^{\star}\in\mathbb{S}}\max_{x}\{\|x-x^{\star}\|:f(x)\leq f(x^0)\}\leq R_0$ and  $\mathbb{S}$ denotes
the optimal solution set. If $f$ is also strongly convex with parameter $\mu >0$, then one has, for  $k=n,2n,3n,\ldots$,
\[
\left( f(x^{k}) \right)-f^\star \leq \left(1-   \frac{\mu}{2\ell_{\max}(1+nL^2/\ell^2_{\max})}\right)^{k/n}\left( f(x^{0}) -f^\star\right).
\]
\end{theorem}

For easy reference, we recall the interpolation conditions from Theorem \ref{T1} in the case that $\mu = 0$:
  The set $\{x^i,\nabla f(x^i),f(x^i)\}$ for $i\in\{0,1,\cdots,N,\star\}$ is $\mathcal{F}_{0,L}-$ interpolable if and only if
  \begin{align}\label{int_glin}
    f(x^i)\geq f(x^j)+\langle\nabla f(x^j),x^i-x^j\rangle+\tfrac{1}{2L}\left\|\nabla f(x^i)-\nabla f(x^j)\right\|^2, \ \ \ \forall i,j\in\{0,1,\cdots,N,\star\}.
  \end{align}
Using these conditions, we may formulate the worst-case convergence rate as performance estimation problem.
\begin{align}\label{PK}
\nonumber   \max & \  f(x^N) -f(x^\star)\\
 \nonumber \st   &  \  \{(x^i; \nabla f(x^i); f(x^i))\} \ \textrm{satisfy \eqref{int_glin} for $i\in\{0,1,\cdots,N, \star\}$ w.r.t.}\ L \\
    & \|x^0-x^\star\|^2\leq\Delta\\
  \nonumber   & \ x^{k} \ \; (k\in\{1,2,\cdots,N\}) \text{ are generated using Algorithm \ref{Alg-cyc}}\\
\nonumber &  x^0\in\mathbb{R}^n, \ \nabla f(x^\star)=0.
\end{align}
Problem \eqref{PK} can be formulated as a semidefinite programming problem, and this is precisely what was done by Kamri et al.\  \cite{kamri2022worst}.

Since the univariate function $t \mapsto f(x^k + te_i)$ is convex and $\ell_i$-smooth,
it follows from \eqref{int_glin} that,
 for every two consecutive points $x^k$ and $x^{k+1}$ generated by
 Algorithm \ref{Alg-cyc}, the following inequalities hold if $i = k \pmod n +1$:
\begin{align}\label{l_ismooth}
  &f(x^k)\geq f(x^{k+1})+\nabla f(x^{k+1})_i(x^k_i-x^{k+1}_i)+\tfrac{1}{2\ell_i}(\nabla f(x^{k})_i-\nabla f(x^{k+1})_i)^2 \\
  &\nonumber f(x^{k+1})\geq f(x^{k})+\nabla f(x^{k})_i(x^{k+1}_i-x^{k}_i)+\tfrac{1}{2\ell_i}(\nabla f(x^{k+1})_i-\nabla f(x^{k})_i)^2.
\end{align}
By adding the above inequalities to \eqref{PK} one can get a better upper bound for the worst-case convergence rate, i.e.
\begin{align}\label{PK_ex}
\nonumber   \max & \  f(x^N) -f(x^\star)\\
 \nonumber \st   &  \  \{(x^i; \nabla f(x^i); f(x^i))\} \ \textrm{satisfy \eqref{int_glin} for $i\in\{1,\cdots,N, \star\}$ w.r.t.}\ L \\
  \nonumber    &  \  \{(x^i; \nabla f(x^i); f(x^i))\} \ \textrm{satisfy \eqref{l_ismooth} for $i\in\{1,\cdots,N, \star\}$ w.r.t.}\ \{\ell_1,\cdots,\ell_n\} \\
    & \|x^0-x^\star\|^2\leq\Delta\\
  \nonumber   & \ x^{k} \ \; (k\in\{1,\cdots,N\}) \text{ are generated using Algorithm \ref{Alg-cyc}}\\
\nonumber &  x^0\in\mathbb{R}^n, \ \nabla f(x^\star)=0.
\end{align}
In order to obtain an SDP relaxation to \eqref{PK_ex}, we proceed in the same way as Kamri et al.\ \cite{kamri2022worst}.
We view \eqref{PK_ex} as a quadratically constrained quadratic program (QCQP) in variables corresponding to the unknowns
\[
x^k_i, \; \frac{\partial f(x^k_i)}{\partial x_i}, \; f(x^k) \quad k\in\{0,1,2,\cdots,N\}, \; i \in \{1,\ldots,n\},
\]
Next we use the following relations to eliminate variables:
\[
x_i^{k+1} - x_i^k = \left\{
\begin{array}{ll}
   -t_k\frac{\partial f(x^k_i)}{\partial x_i} & \mbox{if $i = k \pmod n +1$} \\
  0 & else
\end{array}\right.
\]
which hold for all
$k\in\{0,1,2,\cdots,N-1\}$, and $i \in \{1,\ldots,n\}$.
Subsequently we form the standard Shor SDP relaxation (see e.g.\ \cite{Shor}) of the resulting QCQP. Note that this is different
 to the approach we followed for randomized coordinate descent.
In particular, the size of the SDP relaxation now depends on $n$, which was not the case before.
This also limits the parameter values for which we may solve the SDP relaxations.

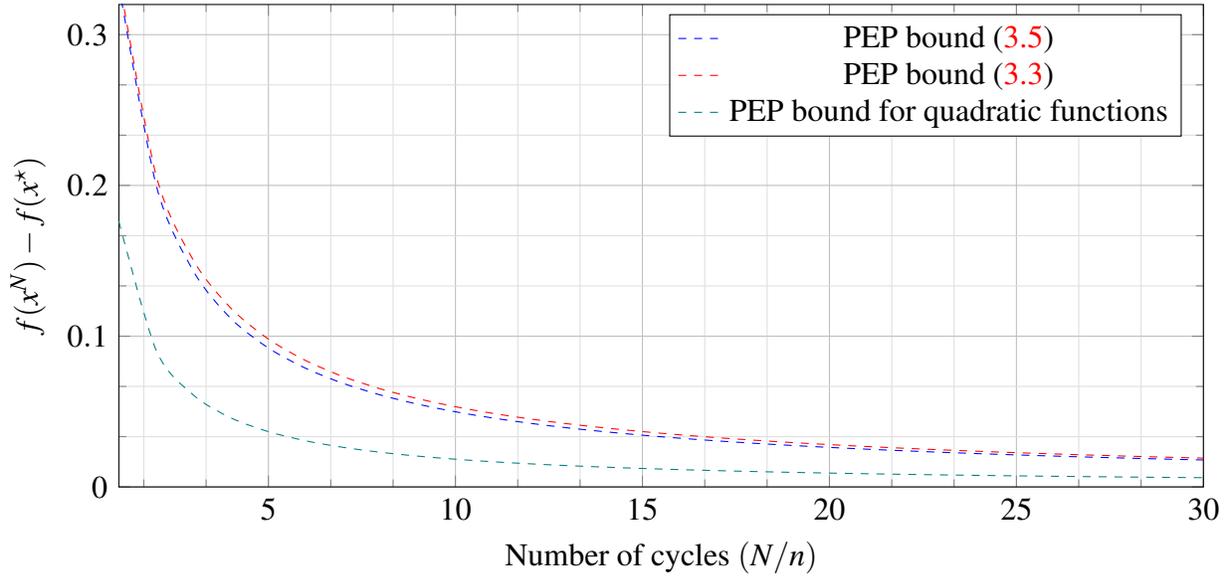
\begin{figure}[H]
\begin{tikzpicture}
\begin{axis}[
    xmin = 1, xmax = 30,
    ymin = 0, ymax = 0.32,
    xtick distance = 5,
    ytick distance = .1,
    grid = both,
    minor tick num = 2,
    major grid style = {lightgray},
    minor grid style = {lightgray!50},
    width = \textwidth,
    height = 0.5\textwidth,
    xlabel = {Number of cycles $(N/n)$},
    ylabel = {$f(x^N) -f(x^\star)$},]

\addplot[
    smooth,
    thin,
    blue,
    dashed
] table[col sep=comma]  {figcoexL2.dat};
\addplot[
    smooth,
    thin,
    red,
    dashed
] table[col sep=comma]  {figcoL2.dat};
\addplot[
    smooth,
    thin,
    teal,
    dashed
] table[col sep=comma]  {figcoquL2.dat};

 \legend{
    PEP bound \eqref{PK_ex},
    PEP bound \eqref{PK},
    PEP bound for quadratic functions,
}

\end{axis}

\end{tikzpicture}
\caption{Convergence rate Algorithm \ref{Alg-cyc} computed by performance estimation problem \eqref{PK_ex} (blue), the bound given by \eqref{PK} (red) and the bound for quadratic functions (green)
for $n=2, L=2, \ell_1=1, \ell_2=1, t=0.5, \Delta=1$.}
\label{fig5}
\end{figure}

In Figure \ref{fig5} we compare the SDP upper bounds from \eqref{PK} and \eqref{PK_ex} for  various parameter values.

The figure shows that the bound can be improved slightly by adding the set of constraints \eqref{l_ismooth}
 to the model provided by Kamri et al.\ \cite{kamri2022worst}. Moreover, we add the constraint which correspond to the quadratic functions
 \eqref{qua_cons} to the model \eqref{PK_ex} which provides us with a better bound  for quadratic functions.
 The computed values are much better that the theoretical bound \eqref{Wrigt2}, to the extent that we do not include this bound in the plot. Indeed,
 Kamri et al.\ \cite{kamri2022worst} already mentioned in their paper that the computed values for their model are much better than the theoretical bound \eqref{Wrigt2}.

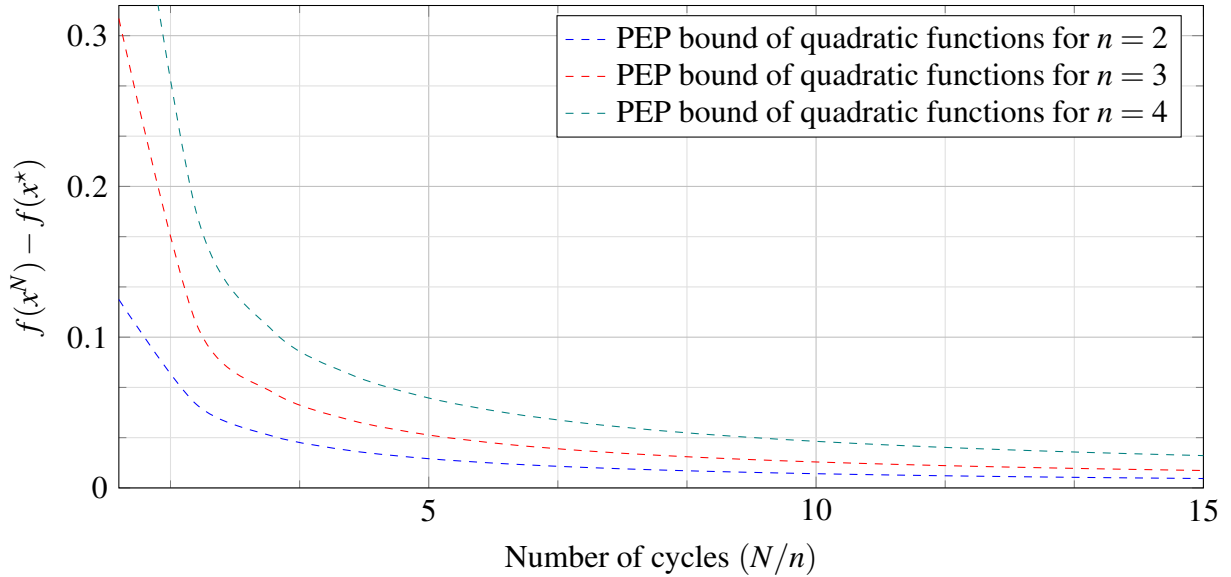
\begin{figure}[H]
\begin{tikzpicture}
\begin{axis}[
    xmin = 1, xmax = 15,
    ymin = 0, ymax = 0.32,
    xtick distance = 5,
    ytick distance = .1,
    grid = both,
    minor tick num = 2,
    major grid style = {lightgray},
    minor grid style = {lightgray!50},
    width = \textwidth,
    height = 0.5\textwidth,
    xlabel = {Number of cycles $(N/n)$},
    ylabel = {$f(x^N) -f(x^\star)$},]

\addplot[
    smooth,
    thin,
    blue,
    dashed
] table[col sep=comma]  {figcop2.dat};
\addplot[
    smooth,
    thin,
    red,
    dashed
] table[col sep=comma]  {figcop3.dat};
\addplot[
    smooth,
    thin,
    teal,
    dashed
] table[col sep=comma]  {figcop4.dat};

 \legend{
    PEP bound of quadratic functions for $n=2$,
    PEP bound of quadratic functions for $n=3$,
    PEP bound of quadratic functions for $n=4$,
}

\end{axis}

\end{tikzpicture}
\caption{Convergence rate for Algorithm \ref{Alg-cyc} computed by performance estimation problem of quadratic functions for $\ell_i=1\ i\in\{1,\cdots,n\}, L=\sum_{i=1}^{n}\ell_i, t=1, \Delta=1$ and different $n$.}
\label{fig6}
\end{figure}

Note that our discussion for coordinate-wise cyclic coordinate descent in this section could be extended to block-wise cyclic coordinate descent in a similar way as was done by Kamri et al.\ \cite{kamri2022worst}.

\subsection{Relation to the Gauss-Seidel method}
The minimization of the convex  quadratic function in \eqref{QP} is  equivalent to the solution of the linear system $Ax = b$. Here, we may assume
w.l.o.g.\ that $A$ has a positive diagonal.
Cyclic coordinate decent for problem \eqref{QP} is closely related to the iterative Gauss-Seidel method for solving this linear system.
For this reason, cyclic coordinate descent is sometimes also referred to as {\em nonlinear Gauss-Seidel}. It is therefore an interesting question whether
the SDP performance estimation framework yields any new insights on the performance of the Gauss-Seidel method.

Denoting $A = (a_{ij})$, the iterative Gauss-Seidel method may be described as follows.
\begin{algorithm}
\caption{Gauss-Seidel method}
\begin{algorithmic}
\State Set $N$  and pick $x^0\in\mathbb{R}^n$.
\State For $k=0, 1, \ldots, N-1$ perform the following:\\
\begin{enumerate}
\item[]
$x_i^{k+1}= \frac{1}{a_{ii}}\left(b_i - \sum_{j=1}^{i-1} a_{ij}x^{k+1}_j- \sum_{j=i+1}^{n} a_{ij}x^{k}_j  \right)$ \quad\quad ($i = 1,\ldots,n$).
\end{enumerate}
\end{algorithmic}
\label{Gauss-Seidel}
\end{algorithm}

This is exactly cyclic coordinate descent with unit step lengths if the gradient at a point $x$ is replaced by $D^{-1}\nabla f(x)$, where $f(x) = \frac{1}{2}x^\top Ax - b^\top x$ as before, and
$D$ is the diagonal matrix with the same diagonal entries as $A$.
To see this, recall that the Fr\'{e}chet derivative of a differentiable function $f:\mathbb{R}^n \rightarrow \mathbb{R}$
  at a point $x \in \mathbb{R}^n$ is the unique linear operator, say $D_f(x):\mathbb{R}^n \rightarrow \mathbb{R}$, such that
\[
\lim_{\| {h}\| \rightarrow 0} \frac{f({x}+{h}) - f({x}) - D_f(x){h} }{\|{h}\|} = 0.
\]
Once an inner product on $\mathbb{R}^n$ is fixed, say $\langle \cdot,\cdot\rangle$, one may, by the Riesz representation theorem, express $D_f(x){h} = \langle g(x),h\rangle$,
where $g(x)$ is called the gradient vector of $f$ at $x$ with respect to $\langle \cdot,\cdot\rangle$. In particular, if $\langle \cdot,\cdot\rangle$
is the Euclidean dot product, then $g(x) = \nabla f(x)$. If one changes to the inner product $\langle \cdot,\cdot\rangle_D$ defined by
\begin{equation}
\label{eq:inner2}
\langle u,v\rangle_D = \sum_{i=1}^n {a_{ii}}u_iv_i \quad\quad (u,v \in \mathbb{R}^n),
\end{equation}
then the gradient vector at $x$ becomes $D^{-1}\nabla f(x)$, by the uniqueness of the Fr\'{e}chet derivative.

It was shown in \cite{de2020worst} that the interpolation condition in Theorem \ref{T1} holds for any reference inner product  $\langle \cdot,\cdot\rangle$,
provided that the gradient vector is interpreted accordingly.

In other words, the following SDP performance estimation problem gives a bound on the worst-case performance of the Gauss-Sidel method after $N$ iterations,
 when $A$ is a
symmetric positive semidefinite matrix with a positive diagonal.

\begin{align}\label{PK_ex2}
\nonumber   \max & \  f(x^N) -f(x^\star)\\
 \nonumber \st   &  \  \{(x^i; \nabla f(x^i); f(x^i))\} \ \textrm{satisfy \eqref{int_glin} and \eqref{qua_cons} for $i\in\{0,1,\cdots,N, \star\}$ w.r.t.\ $L = \lambda_{\max}(D^{-1}A)$} \\
  \nonumber    &  \  \{(x^i; \nabla f(x^i); f(x^i))\} \ \textrm{satisfy \eqref{l_ismooth} for $i\in\{0,1,\cdots,N, \star\}$ w.r.t.\   $\ell_1= \ldots = \ell_n = 1$}\\
    & \|x^0-x^\star\|^2\leq\Delta\\
  \nonumber   & \ x^{k} \ (k\in\{1,2,\cdots,N\}) \text{ is generated using Algorithm \ref{Alg-cyc}}\\
\nonumber &  x^0\in\mathbb{R}^n, \ \nabla f(x^\star)=0,
\end{align}
where the inner product is now understood to be the one in \eqref{eq:inner2}, and the norm the induced norm for this inner product,
and $\lambda_{\max}(D^{-1}A)$ denotes the largest eigenvalue of $D^{-1}A$. (Note that the eigenvalues of $D^{-1}A$ are real.)
Importantly, the reference inner product is not visible in the SDP performance estimation problem reformulation of \eqref{PK_ex2},
since only a Gram matrix for this inner product appears.
It is therefore equally valid, for any inner product, provided that the inner product and norm are interpreted accordingly.
Of course, the Lipschitz constants like \eqref{Co.Lip} depend on the norm as well. It is easy to verify that, for the inner product \eqref{eq:inner2}, and
$f(x) = \frac{1}{2}x^\top Ax - b^\top x$, one has $\ell_1= \ldots = \ell_n = 1$ and $L = \lambda_{\max}(D^{-1}A)$ as is used in \eqref{PK_ex2}.

In summary, we have shown the following.

\begin{theorem}
Consider a solvable system of linear equations $Ax = b$ where $A$ is a symmetric positive semidefinite matrix with positive diagonal, and let $x^\star$ denote a solution.
Letting $f(x) = \frac{1}{2}x^\top Ax - b^\top x$, after $N$ iterations of the Gauss-Seidel method, an
upper bound on $f(x^N) - f(x^\star)$ is given by the optimal value of the SDP problem \eqref{PK_ex2}, provided
that the starting point $x^0$ satisfies $\|x^0 - x^\star\| \le \Delta$ for a given $\Delta$, where the norm is the induced norm of the inner product \eqref{eq:inner2}.
\end{theorem}

The Gauss-Seidel method is known to be convergent when $A$
is symmetric positive-definite, e.g.\ \cite[Theorem 10.1.2]{Golub1}, or strictly or irreducibly diagonally dominant, e.g.\ \cite{bagnara1995unified}. The case when $A$ is only positive semidefinite
(with positive diagonal) seems to be less well-understood, and our approach sheds more light on this case. In particular,
 numerical results of the type shown  in Figure \ref{fig6} apply here.

\section{Conclusion}
We have studied SDP performance estimation approaches to analyse randomized and cyclic coordinate descent,
 thereby complementing recent results in \cite{kamri2022worst}.
 For randomized coordinate descent, we have given the first known SDP performance estimation bound.
 For cyclic coordinate descent, we were able to improve slightly on the numerical values given in \cite{kamri2022worst}, and we also
 discussed the link with the Gauss-Seidel method in the case of convex quadratic functions.
Of course, to obtain new rates of convergence in general, it is necessary to solve the  SDP performance estimation problems analytically, as opposed to numerically, but we have been unable to obtain analytic solutions for the SDP problems presented in this paper.

\vskip 6mm
\noindent{\bf Acknowledgments}

\noindent   This work was supported by the Dutch Scientific Council (NWO)  grant OCENW.GROOT.2019.015, \emph{Optimization for and with Machine Learning (OPTIMAL)}.

\bibliographystyle{unsrt}
\bibliography{Bibliography}


\end{document}